\documentclass[12pt]{article}
\usepackage{graphicx}

\def\ds{\displaystyle}
\begin{document}

\title{ A class of Lorenz-type systems, their factorizations and extensions}
\author{I.~Kunin$^1$ and A.~Runov$^2$}
\maketitle

\footnotetext[1]{
	 Department of Mechanical Engineering, 
	 University of Houston,
	 Houston, TX 77204, USA, e-mail kunin@uh.edu
	 }
\footnotetext[2]{
	 Department of Theoretical Physics, 
	 St.~Petersburg University,
	 Uljanovskaja~1, St.~Petersburg, Petrodvorez, 198504, Russia
	 }

\begin{abstract}
It is well known that the Lorenz system has $Z_2$-symmetry. Using  introducted
in~\cite{aa_duff} topological covering-coloring a new
representation
for the Lorenz system  is obtained.  Deleting coloring leads to the factorized
Lorenz system that is in a sense more fundamental than the original one.
Finally, $Z_n$ extensions define a class of Lorenz-type systems. The approach
admits a natural generalization for regular and chaotic systems with arbitrary
symmetries. 
\end{abstract}

There are more than a thousand publications related to the Lorenz system and
the number is increasing. Some explanations for this popularity are given,
e.g. in~\cite{sparr, jacks}. Here  the Lorenz system is realized in a
non-traditional
form that may motivate a new approach to investigating chaotic systems with
symmetries. The idea of this approach has been demonstrated in~\cite{aa_duff} for the
simple and instructive example of the Duffing oscillator with $Z_2$-symmetry. The
basic components are: realization of the symmetry by topological
covering-coloring, factorization by deleting  the coloring, and the following
$Z_n$-extension. As a result a class of Lorenz-type systems is
defined. 
An open and intriguing problem: comparison of chaos for
elements of the class. 

\begin{figure}[p]
\centerline{
\includegraphics[width=10cm]{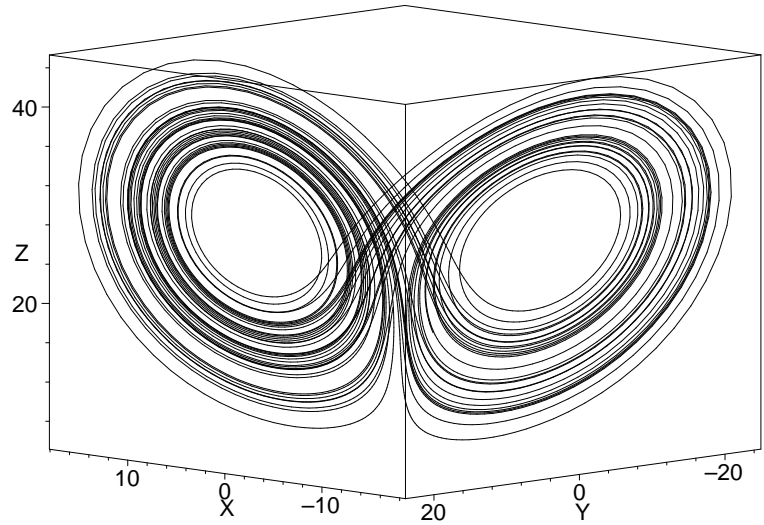}
}
\caption{Standard Lorenz system.}
\label{stdLor}
\end{figure}

\begin{figure}[p]
\centerline{
\includegraphics[width=10cm]{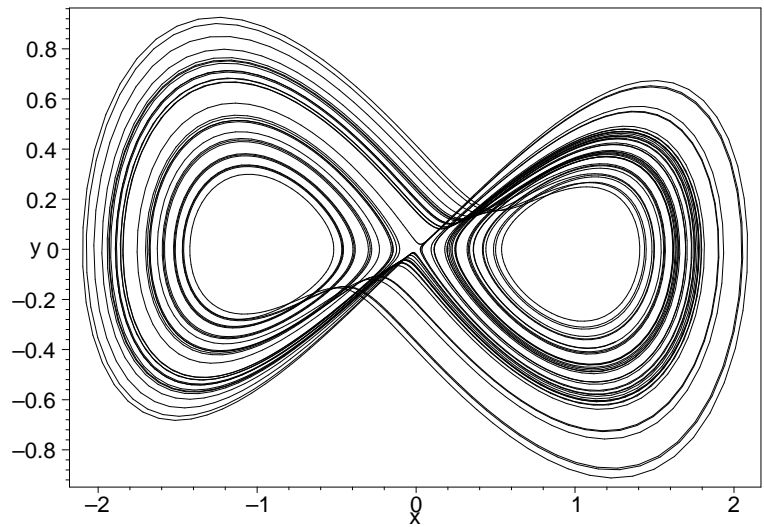}
}
\caption{Rescaled Lorenz system (L$_2$).}
\label{rLor}
\end{figure}

\begin{figure}[p]
\centerline{
\includegraphics[width=10cm]{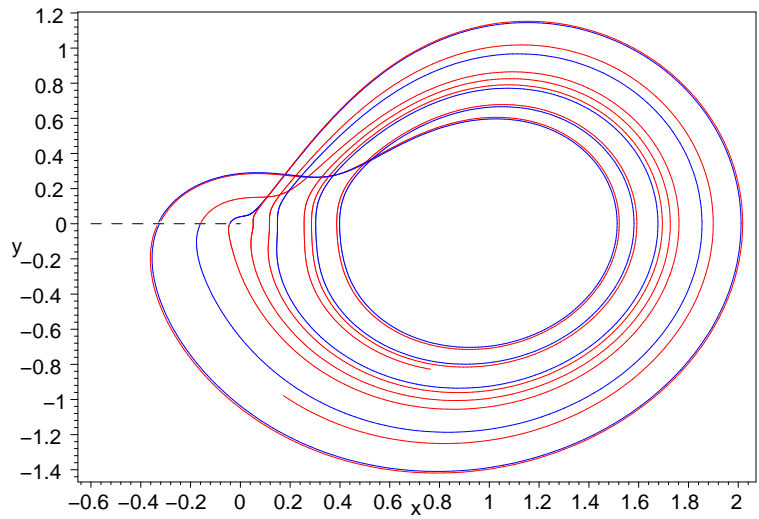}
}
\caption{Covered (colored) L$_2$ system, $x_1$ --- $y_1$ projection, (short piece of trajectory).}
\label{cLor2D}
\end{figure}

\begin{figure}[p]
\centerline{
\includegraphics[width=10cm]{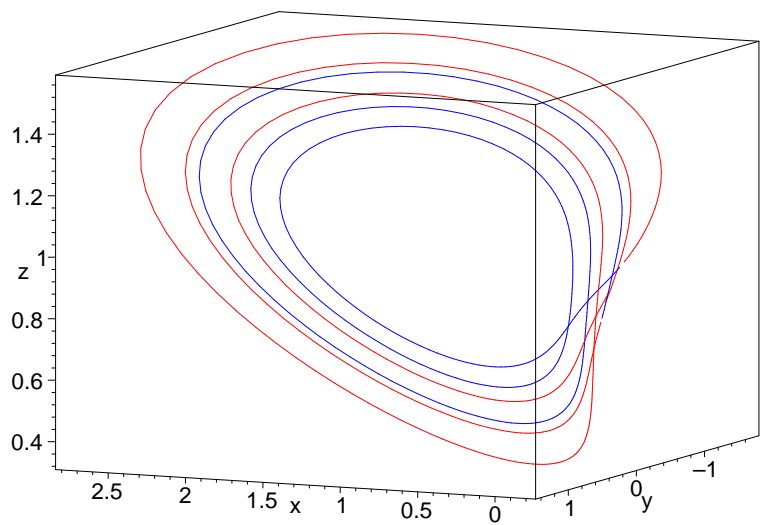}
}
\caption{Covered (colored) L$_2$ system, 3-D projection,(short piece of
trajectory).}
\label{cLor3D}
\end{figure}

Let us start with the standard form of the Lorenz system
\begin{equation} \label{SL}
 \left\{
 \begin{array}{l}
 \dot X=\sigma(Y-X)\\
 \dot Y=(r-Z)X-Y\\
 \dot Z=-bZ+XY
 \end{array}
 \right.
\end{equation}
and the corresponding canonical strange attractor (Fig.~\ref{stdLor}). It is
convenient to
use the following normalized form L$_2$
\begin{equation} \label{RL}
 \left\{
 \begin{array}{l}
 \dot x=y\\
 \dot y=((1-\gamma) x^2-1+z)x - \mu y\\
 \dot z=\beta(\gamma x^2-z)
 \end{array}
 \right.
\end{equation}
related to~(\ref{SL}) by  (for $r>1$)
\begin{equation}
 \begin{array}{l}
 t=\ds\sqrt{(r-1)\sigma} t_L \\
 x=\ds\frac1{\sqrt{(r-1)b}} X \\
 y=\ds\frac1{r-1}\sqrt{\frac\sigma b} (Y-X) \\
 z=\ds\frac1{r-1} \left(Z-\frac1{2\sigma}X^2 \right) \\
 \mu=\ds\frac{1+\sigma}{\sqrt{(r-1)\sigma}}, \quad
 \beta=\ds\frac b{\sqrt{(r-1)\sigma}}, \quad
 \gamma=\ds 1-\frac b{2\sigma}
 \end{array}
\end{equation}
The system L$_2$ has 3 fixed points: $(0,0,0)$, $(\pm1,0,\gamma)$
and $Z_2$-symmetry: $(x,y) \to (-x,-y)$. Corresponding attractor is shown on
Fig.~\ref{rLor}. Here we have an analogy with the Duffing
system~\cite{aa_duff}
motivating the similar transformation of topological covering 
\begin{equation}
 \begin{array}{l}
  x_1=x_1(x,y)=\displaystyle\frac{x^2-y^2}r,     \\
  y_1=y_1(x,y)=\displaystyle\frac{2xy}r.
 \end{array} \qquad r=\sqrt{x^2+y^2}
\end{equation}
The corresponding colored trajectories are shown in
Figs. \ref{cLor2D},~\ref{cLor3D}. 
	
The meaning of the transformation becomes more transparent in polar
coordinates: it is $(r,\phi) \to (r,2\phi)$. The symmetry of the
system L$_2$ is equivalent to $f(r,\phi)=f(r,\phi+\pi)$, where f is an arbitrary
function on the states of L$_2$. Correspondingly, in new coordinates
$f(r,\phi)=f(r,\phi+2\pi)$ , i.~e. one deals with two identical planes marked
by colors. The
equations for this representation of the original Lorenz system are
uniquely defined by the indicated above transformations. Their
explicit form  can be easily obtained  but is far from being elegant
(really it is  cumbersome) and is not given here.

\begin{figure}[p]
\centerline{
\includegraphics[width=10cm]{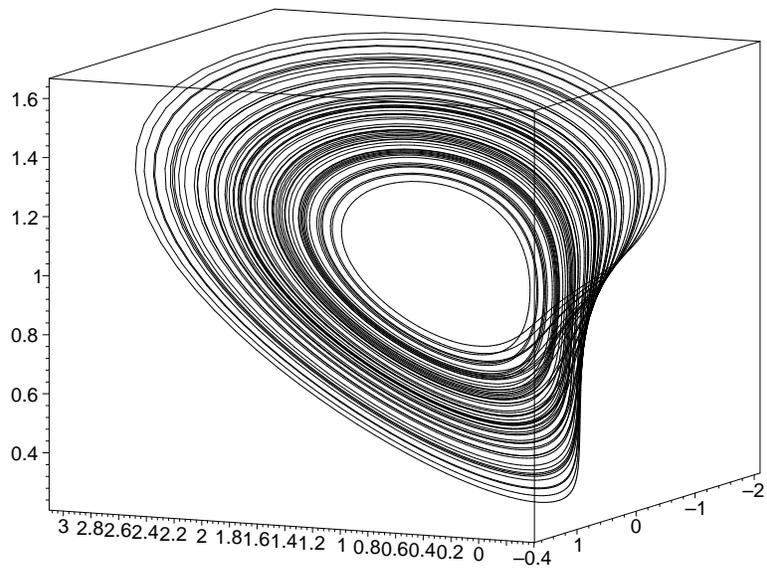}
}
\caption{Factorized Lorenz attractor (L$_1$).}
\label{FLor}
\end{figure}

\begin{figure}[p]
\centerline{
\includegraphics[width=10cm]{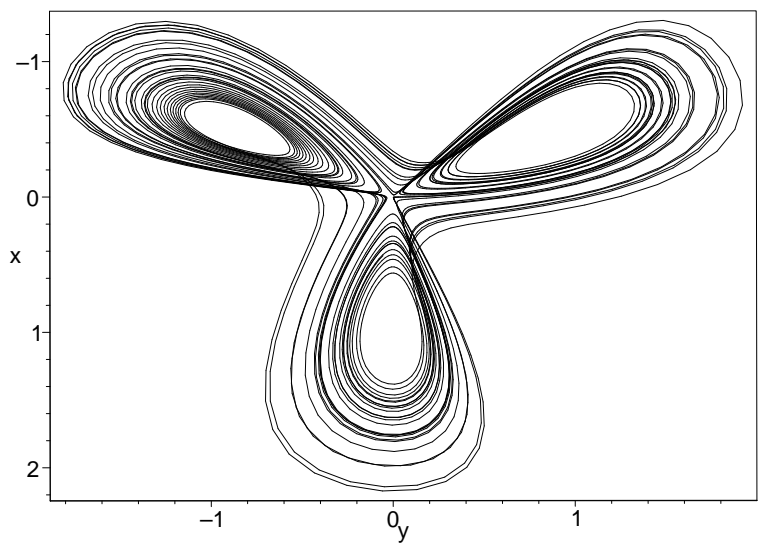}
}
\caption{Attractor of L$_3$.}
\label{L3}
\end{figure}

The erasing of the colors is equivalent to gluing together two underlying
planes, i.~e. defining a new, more simple and fundamental factorized Lorenz
system L$_1$ (Fig.~\ref{FLor}). As above the explicit equations can be
written, but  they
are cumbersome, and we do not present them here. The  system L$_1$ carries all
essential information about L$_2$ except the symmetry. Thus the Lorenz system may
be naturally called a Z$_2$-extension of the L$_1$ system. This motivates to
define L$_n$ as Z$_n$-extension of L$_1$, i.~e. the similar transformation to the
corresponding system with Z$_n$-symmetry. The example for L$_3$ is given in
Fig.~\ref{L3}.
Thus we have defined  the class of Lorenz-type systems. Any two elements of
the class are related by compositions of appropriate
factorizations and extensions. It is natural to call these
elements a class of 
Z-equivalent systems (with a canonical representative L$_1$).

The approach admits a generalization to different families of discrete and
continuous symmetries as well as wider classes of transformations.
The scheme of such generalizatins will be given elsewhere.

\end{document}